\documentclass[12pt]{article}

\newlength{\stefan}
\usepackage{amsmath, amsthm, amsfonts, makeidx}
\usepackage{multind}
\DeclareMathSymbol{\subsetneq}{\mathord}{AMSb}{"26}

\newtheorem{lemma}{Lemma}[section]

\newtheorem{theorem}[lemma]{Theorem}

\newtheorem{corollary}[lemma]{Corollary}
\theoremstyle{definition}

\newcommand{\lp}{\longrightarrow}

\newcommand{\mb}{\mathbb}

\newcommand{\F}{\mathcal{F}}

\newcommand{\X}{\mathcal{X}}

\newcommand{\C}{\mb{C}}

\newcommand{\N}{\mb{N}}

\newcommand{\Aut}{\mathit{Aut}}

\newcommand{\LND}{\operatorname{LND}}

\newcommand{\spec}{\operatorname{spec}}

\newcommand{\ttr}{\tilde{r}}
\newcommand{\tts}{\tilde{s}}

\title{On the methods to construct UFD counterexamples to a cancellation problem}
\author{Stefan Maubach\footnote{Funded by Veni-grant of council for the
physical sciences, Netherlands Organisation for scientific research (NWO)}\\ \ \\
\small
Radboud University Nijmegen\\\small Toernooiveld 1, The Netherlands\\ \small s.maubach@math.ru.nl}

\begin{document}

\maketitle

\begin{abstract}
In a previous paper  \cite{FM06} the author together with prof. dr. Finston constructed a class of UFDs 
$A_{n,m}$ where $n,m\in \N^*$. These rings are all stably equivalent 
($A_{n,m}[T]\cong  A_{p,q}[T]$ for all $n,m,p,q$) but are only isomorphic themselves if $(n,m)=(p,q)$. 
These examples are the first UFD examples over a characteristically closed field satisfying this behavior. 
In this paper, we describe the methods used in this article, and show that they are very general, enabling the 
reader to construct many more such examples, based on the same principles.
\end{abstract}

\section{Introduction}

This paper zooms in on what is essential in the example in the paper \cite{FM06}.
Let us repeat a typical example of this paper: (we write $R^{[1]}$ for a polynomial ring in one variable over $R$.)

Define $R:\C[x,y,z]:=\C[X,Y,Z]/(X^2+Y^3+Z^7)$, and let $A_{n,m}:=R[u,v]=R[U,V]/(x^mU-y^nV-1)$ where $n,m$ are positive integers.
Now it is shown in \cite{FM06} that $A_{n,m}^{[1]}\cong A_{n',m'}^{[1]}$  for any positive integers $n,m,n',m'$, while
$A_{n,m}\cong A_{n',m'}$ implies that $(n,m)=(n',m')$. This is a UFD-counterexample to the so-called {\em generalized cancellation problem}, 
which states: does $R^{[1]}\cong S^{[1]}$ imply that $R\cong S$? The mentioned example is the ``best worst'' example yet, being the ``nicest'' rings  $R$ and $S$ 
for which the generalized cancellation problem does not hold. The big conjecture at the moment is what nowadays is called ``the'' cancellation problem: 
the case that $S=\C^{n]}$. I.e. does $R^{[1]}\cong \C^{[n+1]}$ imply $R\cong \C^{[n]}$? (This problem is still open for $n\geq 4$.)

However, it seems like in this type of counterexample to the generalized cancellation problem, the ring $R$ can be chosen much more freely. 
For a ring $R$ and elements $r,s\in R$, write $A_{r,s}:=R[U,V]/(rU-sV-1)$.
So we are looking for a ring $R$ and elements $r,s,r',s'$ in $R$ such that
(1) $A_{r,s}\not \cong A_{r',s'}$, while $A_{r,s}^{[1]}\cong A_{r's'}^{[1]}$, (2) $A_{r,s}$ and $A_{r's'}$ are $\C$-algebra UFDs of dimension 3.\\
It is not our goal to classify which rings $R$ have elements $r,s,r',s'$ having the above properties, but we want to discuss properties
that enable us to give examples. These properties are mainly for the part of showing that $A_{r,s}$ is not isomorphic to $A_{r',s'}$, except \ref{SS6}.

\subsection{Notations}

{\bf Notations:} If $R$ is a ring, then $R^{[n]}$ denotes the polynomial ring in $n$ variables over $R$.
We will use the letter $k$ for a field of characteristic zero, and $K$ for a fixed algebraic closure.
When $X,Y,\ldots$ are variables in a polynomial ring of rational function field, denote by  $\partial_X, \partial_Y,\ldots$ the derivative with respect to $X,Y,\ldots$.
Very often, we will use small caps $x,y,z,\ldots$ for residue classes of  $X,Y,Z,\ldots$ modulo some ideal.

\section{Useful properties of the rings $R$ and $A_{r,s}$}

\subsection{$R$ must be a UFD, and $A_{r,s}$ must be a UFD.}
\label{SS1}

It is not true that $R$ {\em must} be a UFD to make $A_{r,s}$ into a UFD. For example, if $R_{p,q}:=\C[X,Y,Z]/(X^pY-Z^q)$ and $A_{p,q,m,n}=R_{p,q}[U,V]/(x^mU-y^nV-1)$
then one can show that $A_{p,q,m,n}\cong \C[X,Z,V,X^{-1}]$ for any choice of $p,q,m,n\in \N^*$, which is a UFD.(Proofsketch: $A_{p,q,m,n}$ can be seen as a subring of 
$\C[X,Z,V,X^{-1}]$ where $Y=Z^qX^{-p}$ and $U=(Y^nV+1)X^{-m}$. Define $\tilde{Y}:=X^{p-1}Y, \tilde{U}:=X^{m-1}U$. 
If $q>n$ consider $\tilde{U}-\tilde{Y}VX^{q-n}$, if $q\leq n$ then $X^{-1}=\X^{q-n}\tilde{U}-\tilde{Y}V$. 

Even though $R$ does not need to be a UFD, we require it as computations are much easier (it might be dropped, though).
In order to prove that a ring is a UFD, it is sometimes necessary to compute the class  group (see \cite{Fossum}). The class field
group tells one ``how far'' a ring is from being a UFD, as being a UFD is equivalent to the class group being trivial, for integrally closed noetherian rings. It is not always an easy task to do that,
however. We will quote a few useful tools:

\begin{theorem}
\label{Fossum} (Corollary 10.3 of \cite{Fossum}) Let $A=A_0+A_1+\ldots$ be a
graded noetherian Krull domain such that $A_0$ is a field. Let $\mathfrak{m}
=A_1+A_2+\ldots$. Then $Cl(A)\cong Cl(A_{\mathfrak{m}})$, where $Cl$ is the
class group.
\end{theorem}

\begin{theorem}\label{HO}
\label{local}(\cite{HO74}) A local noetherian ring $(A,\mathfrak{m})$ with
characteristic $A/\mathfrak{m}=0$ and an isolated singularity is a UFD if
its depth is $\geq 3$ and the embedding codimension is $\leq dim(A)-3$.
\end{theorem}

The latter two theorems can be used to show that the hypersurface $X_1^{d_1}+X_2^{d_2}+\ldots+X_n^{d_n}$ is factorial if $n\geq 5$ and any $d_i\in \N^*$
(see for example \cite{FM07} for a proof).
However, theorem \ref{HO} is not that useful here, if one wants to have a 2-dimensional UFD.

One of the more straightforward tools is

\begin{theorem}(Nagata)
Let $A$ be a domain, and let $x\in A$ be a prime element. If $A[x^{-1}]$ is a UFD, then $A$ is a UFD.
\end{theorem}

This is especially useful in showing that $A_{r,s}$ is a UFD, depending on what $r$ and $s$ are.

\begin{lemma}
Let $r$ or $s$ be a prime element in $R$, assume $R$ is a noetherian UFD, and assume $r$ and $s$ share no common factor. Then $A_{r,s}$ is a UFD.
\end{lemma}

\begin{proof}
Write $r=r_1r_2\ldots r_k$ where the $r_i$ are irreducible (which can be done since $R$ is noetherian) and prime (which follows since $R$ is a UFD).
We will proceed by induction to $k$. 
If $k=0$ then $r$ is invertible and $A_{r,s}\cong R[V]$. 

Now $r_k$ is prime in $A_{r,s}$, since $A_{r,s}/(r_k)\cong R[U,V]/(r_k, -sV-1) =(R/r_k)[1/({s\mod r_k})]$ which is a domain.
$A_{r,s}[r_k^{-1}]=R[r_k^{-1}][U,V]/(rU-sV-1)$, which is a UFD by induction (as $r\in R[r_k^{-1}]$ has fewer irreducible factors) and Nagata's theorem. 
\end{proof}

\subsection{$R^*=A^*$}

\label{SS2}

This also implies that $r$ and $s$ do not share a common factor other than a unit, as this common factor will become invertible in $A_{r,s}$.

\subsection{$R$ is rigid, $ML(A_{r,s})=R$}
\label{SS3}

$R$ being rigid is defined as $\LND(R)=\{0\}$, i.e. there are no nontrivial $G_a$-actions on the variety associated to $R$.
An equivalent definition is that the Makar-Limanov invariant is maximal, i.e.$ML(R)=R$.
This is not a necessary property for a counterexample to generalized cancellation, but it is very useful in making sure that $A_{r,s}$ has few automorphisms.
Interesting to note is that this is the point where we already rule out the possibility for constructing a counterexample to ``the'' cancellation problem,
as having few automorphisms contradicts being isomorphic to $\C^{n}$.  The reason that we require this here, is that
we will want to distinguish $A_{r,s}$ and $A_{r',s'}$ later on by computing their automorphism groups.
Also, this will automatically take care of the next requirement.

In order to make a rigid ring, we bump into a strange phenomenon. It seems like ``almost any'' ring is rigid, but it is in general hard to prove that a ring
is rigid. Note also that, through this difficulty, it is very dangerous to make statements as ``almost any'' ring is rigid, as it is hard to prove any such statement.
On a side note, no examples are known of rigid rings $R$ for which $ML(R^{[n]})\not = ML(R)$, we refer to \cite{Crach, CML05} for comments on this difficult problem (``losing rigidity''). This is connected with the additional requirement that $ML(A_{r,s})=R$: we have an extension $A$ of the rigid ring $R$, and in general,
$ML(A)$ can be anything: equal to $R$, strictly containing $R$ (like being rigid itself), and we even cannot exclude $ML(A)$ being strictly contained in $R$.
Note that, in this case, we do have $ML(A_{r,s})\subseteq R$ as $s\partial_u+r\partial_v\in LND(A_{r,s})$, which has kernel $R$ as can be easily checked.
Here we can view $\partial_u$ (resp. $\partial_v$) as the restrictions to $A_{r,s}$ of the partial derivatives with respect to $U$ (resp. $V$) on  
$Q(R)[U]$ (resp. $Q(R)[V]$) where $Q(R)$ denotes the quotient field of $R$.

There are a few ways of constructing and proving that a ring is rigid. A very useful lemma is the following (lemma 2.2 in \cite{FM07}):

\begin{lemma}
\label{L1} Let $D$ be a nonzero locally nilpotent derivation on a domain $A$
containing $\mathbb{Q}$. Then $A$\ embeds into $K[S]$ where $K$ is some
algebraically closed field of characteristic zero, in such a way that $%
D=\partial _{S}$ on $K[S]$.
\end{lemma}

For example: one has a domain $R:=\C^{[n]}/(F)$ where $F\in \C^{[n]}$. If there exists some nontrivial $D\in \LND(R)$, then we can see the elements and also variables
 of $R$ as elements in $K[S]$. So, $F=0$, but also $0=\partial_S(F)=\sum (\partial_S X_i(S)) \frac{\partial F}{\partial X_i}$. These two equations
 can yield that each $X_i(S)$ is constant in $S$. If that is the case, then $D$ is the zero map, and one has a contradiction. This is exploited in both
\cite{FM06} and \cite{FM07}, using (an extension of) Mason's Theorem.

Incidentally, one can also use this method to construct rings with a restricted supply of LNDs. See \cite{FM07} and \cite{MauPre}.

\subsection{$R$ must be a characteristic subring of $A_{r,s}$}
\label{SS4}

A characteristic subring is a subring which stays invariant under all automorphisms. If $ML(A_{r,s})=R$, then $A_{r,s}$ will have this property:

\begin{lemma}
The Makar-Limanov invariant of a ring $B$ is a characteristic subring of $B$.
\end{lemma}

For a proof, see for example \cite{FM06} lemma 4. This does imply that

\begin{corollary}\label{corr}
Any $\varphi\in Aut_{\C}(A_{r,s})$ satisfies $\varphi(R)=R$.
\end{corollary}

\begin{lemma}
$\LND(A_{r,s})= R E $ where $E=s\partial_u+r\partial_v$.
\end{lemma}

\begin{proof}
Since $ML(A_{r,s})=R$, any $D\in \LND(A_{r,s})$ will satisfy $D(r)=D(s)=0$. Therefore, $0=D(ru-sv-1)$ implies
$rD(u)=sD(v)$. Now here it is handy if one knows $A_{r,s}$ to be a UFD (otherwise the following may still be true, but much more complicated)
as we can conclude that $D(u)=st, D(v)=rt$ for some $t\in A_{r,s}$ (since $r,s$ share no common factor).
So $D=tE$, and now we can use the well-known result that if $fD\in LND(B)$ for some ring $B$, then $D\in LND(B)$ and $D(f)=0$.
This implies $D\in RE$.
\end{proof}

\subsection{The restriction $\F: \Aut_{\C}(A_{r,s}) \lp \Aut_{\C}(R)$ must be surjective}
\label{SS5}
\label{Subs}
Note that this restriction $\F$ exists because of corollary \ref{corr}. What we require here is surjectivity.
This property moves the problem to determining $\Aut_{R}(A_{r,s})$.

\subsection{$(r,s)$ is a height 2 ideal of $R$}

\label{SS6}
We will need in lemma \ref{AutoAR} that $(r,s)\not = R$, which is implied by this requirement, but we mainly need this requirement 
for the following:

\begin{lemma} \label{gelijk}
If $rad(r,s)=rad(r',s')$ then $A_{r,s}^{[1]}\cong A_{r',s'}^{[1]}$.
\end{lemma}

\begin{proof}
Let us write $X_{r,s}$ for the variety associated to $A_{r,s}$. 
We have a $G_a$-action on $A_{r,s}$ (associated to $s\partial_u+r\partial_v$).

The $G_a$-action is locally trivial 
(in fact the basic open subsets $\mathcal{D}_X(r)$ and $\mathcal{D}_X(s)$ cover $X_{r,s}$ and satisfy  
$\mathcal{D}_X(s)=\mathcal{D}_{\spec R} (s) \times \C$,  $\mathcal{D}_X(r)=\mathcal{D}_{\spec R} (r) \times\C$). 
Therefore $X_{r,s}$ is the total space of an algebraic principal $G_a$-bundle over  $\spec(R)\backslash \mathcal{V}$ where $\mathcal{V}$ is the
set of all prime ideals containing $(r,s)$. The same for $X_{r',s'}$.
Now we can take their fiber product over the base:
$X_{r,s}\times_{\spec{(R)}\backslash \mathcal{V}} X_{r',s'}$.
By standard arguments, since $X_{r,s}$ and $X_{r',s'}$ are affine, this is isomorphic to $X_{r,s}\times \C$ as well as $X_{r',s'}\times \C$.
So $A_{r,s}^{[1]}=\mathcal{O}(X_{r,s}\times \C)=\mathcal{O}(X_{r',s'}\times \C)=A_{r',s'}^{[1]}$.
\end{proof}

\section{The $R$- automorphism group of $A_{s,t}$}

If one has $R, A_{r,s}$ satisfying everything in the previous section, then there are some things which come for free.
To be more precise, $\Aut_{R}(A_{r,s})$ can be described, and we can give a simple requirement such that $A_{r,s}\not \cong A_{r',s'}$.

\begin{lemma} \label{05.10.cor4} Let $\varphi\in  Aut_{\C}(A_{r,s})$. Then $\varphi^{-1}E\varphi = \lambda E$ where $\lambda\in R^*$.
\end{lemma}

\begin{proof}
$\varphi^{-1}(LND(A_{r,s}))\varphi=\LND(A_{r,s})$, as can be easily proved since conjugating an LND yields another LND (showing $\subseteq$), and conjugating with $\varphi^{-1}$ gives $\supseteq$.
Therefore, $RE=R(\varphi^{-1}E\varphi)$ and the result follows.
\end{proof}

\begin{lemma}\label{AutoAR}
{$\varphi\in Aut_{R}A_{n,m}$ if and only if $\varphi$ is an $R$-homomorphism
satisfying $\varphi(u,v)=(ts+u,tr+v)=exp(tE)$ for some
$t\in R$. Consequently, $Aut_{R}A_{n,m}\cong<R,+>$ as groups. }
\end{lemma}

\begin{proof}
{ We know by corollary \ref{05.10.cor4} that $\varphi^{-1}(E)\varphi=\lambda
E$ for some $\lambda\in R^*$. Define $(F,G):=(\varphi(u),\varphi(v))$ and denote this by $\varphi(u,v)$. Similarly, $E(F,G):=(E(F),E(G))$. Also,
$\varphi|_R=Id.$. So now
\[%
\begin{array}
[c]{rl}%
(\lambda s,\lambda r)= & \varphi(\lambda s,\lambda r)\\
= & \varphi\lambda E(u,v)\\
= & \varphi(\varphi^{-1}E\varphi)(u,v)\\
= & E(F,G)\\
= & (sF_{u}+rF_{v},sG_{u}+rG_{v})
\end{array}
\]
where the subscript denotes$\ $partial derivative.}

{ Let us consider the first equation, }%
\[
{\lambda s=sF_{u}+rF_{v}.}%
\]
{ Defining $H:=F-\lambda u$, we see that $-sH_{u}=rH_{v}$. By the
following lemma \ref{05.10.diffeq} we see that $H=p\in R$, so }%
\[
{F=p+\lambda u.}%
\]
{ The second equation yields $\lambda r=sG_{u}+rG_{v}$. Defining
$H:=G-\lambda v$, yields $-rH_{v}=sH_{u}$, which by the following
lemma \ref{05.10.diffeq} yields $H=q\in R$ and thus $G=q+\lambda v$.
Now}%

\[%
\begin{tabular}
[c]{lll}%
$0$ & $=$ & ${\varphi(ru-sv-1)}$\\
& $=$ & ${r\varphi(u)-s\varphi(v)-1}$\\
& $=$ & ${rF-sG-1}$\\
& $=$ & ${r(p+\lambda u)-s(q+\lambda v)-1}$\\
& $=$ & ${rp-sq+\lambda(ru-sv)-1}$\\
& $=$ & ${rp-sq+\lambda-1.}$%
\end{tabular}
\
\]

Now due to \ref{SS6}, $1-\lambda=rp-sq$ are in a maximal ideal, hence $\lambda=1$.
Therefore, $rp=sq$, and since $r$ and $s$ share no common factor, and $R$ is a UFD, we get that
$p=st$ and $q=rt$ for some $t\in R$. Thus any automorphism must have the given form.
It is not difficult to check that maps of this form are well-defined
homomorphisms which are automorphisms.
\end{proof}

\begin{lemma}
{ \label{05.10.diffeq} If $H\in A_{r,s}$ such that $-sH_{u}=rH_{v}$,
then $H\in R$. }
\end{lemma}

\begin{proof}
{ We can find polynomials $p_{i}(v)\in$}${R[v]}${
such that $H=\sum_{i=0}^{d}p_{i}u^{i}$ for some $d\in\mathbb{N}$. Requiring
that $r$ does not divide coefficients of $p_i(v)$ if $i\geq 1$ (which we can do as $ru=sv+1$) we force the  $p_i$ to be unique.
The
equation $-y^{n}H_{u}=x^{m}H_{v}$ yields
\[
\sum_{i=0}^{d-1}-(i+1)s p_{i+1}u^{i}=\sum_{i=0}^{d}r p_{i,v}u^{i}%
\]
where }$p_{i,v}\equiv\frac{\partial p_{i}}{\partial v}.$ \ {Substitute
$sv+1$ for $ru$ to obtain a unique representation:}%
\[%
\begin{array}
[c]{rl}%
\sum_{i=0}^{d-1}-(i+1)sp_{i+1}u^{i}= & rp_{0,v}+\sum_{i=0}^{d-1}%
(sv+1)p_{i+1,v}u^{i},
\end{array}
\]
so
\[%
\begin{array}
[c]{rl}%
-sp_{1}= & rp_{0,v}+(sv+1)p_{1,v}%
\end{array}
\]
and%
\[%
\begin{array}
[c]{rl}%
-(i+1)sp_{i+1}= & (sv+1)p_{i+1,v}%
\end{array}
\]
for each $i\geq1$.

{ Let $i\geq1$ and assume that $p_{i+1}$ has degree $k$ with respect to $v$.
Let $\alpha \in R$ be the top coefficient of $p_{i+1}$, seen as a polynomial
in $v$. Then $-(i+1)s\alpha=sk\alpha$, but that gives a contradiction.
So for each $i\geq1:p_{i+1}=0$. This leaves the equation $0=rp_{0,v}$
which means that $p_{0}\in R$. Thus $H=p_{0}u^{0}\in R$. }
\end{proof}

\begin{theorem} \label{ongelijk}Let $R$, $A_{r,s}$, $A_{r',s'}$ satisfy the requirements of the previous section.
Suppose that $A_{r,s}\cong A_{r',s'}$. Then there exists $\varphi\in \text{Aut}_{\C}(R)$ such that 
$\varphi(r)R+\varphi(s)R=r'R+s'R$.
\end{theorem}

\begin{proof}
Let $\sigma:A_{r,s}\lp A_{r',s'}$ be an automorphism. Since $\sigma(ML(A_{r,s})=ML(A_{r',s'})$ we know that $\sigma(R)=R$. 
Since any automorphism of $R$ is the restriction of an automorphism of $A_{r',s'}$ by \ref{Subs} 
(this is exactly the spot where we use this requirement), we can compose $\sigma$ by an appropriate automorphism $\varphi$ of $A_{r',s'}$, 
and can assume that $\Phi:=\sigma\varphi$ is the identity on $R$. Write $\ttr:=\varphi(r), \tts:=\varphi(s)$. 

Now set $K:=Q(R)$, the quotient field of $R$. Identify $K \otimes_R A_{\ttr,\tts}$ with $K[v]$, $K \otimes_R A_{r',s'}$ with $K[v']$,
and note that $\Phi$ can be extended to a $K$-isomorphism $K[v]\lp K[v']$. So we can assume that $\Phi(v)=\alpha v'+\beta$ where $\alpha\in K^*,\beta\in K$.

Of each ring $A_{\ttr,\tts}$ and $A_{r',s'}$ we know the set of locally nilpotent derivations. Let $\LND(A_{\ttr,\tts})=RE$ and $\LND(A_{r',s'})=RE'$, where
$E(u)=\tts,E(v)=\ttr, E'(u')=s', E'(v')=r'$. Since $\Phi^{-1}\LND(A_{r',s'})\Phi=\LND(A_{\ttr,\tts})$, we must have $\Phi^{-1}E'\Phi=\lambda E$ where
$\lambda\in A^*_{\ttr,\tts}=R^*$.

A computation shows that
\[ \lambda \ttr =\lambda E(v)= \Phi^{-1}E'\Phi (v) = \alpha r' \]
and thus $\alpha=\lambda \ttr/r'$.

Now $\alpha V'+\beta \in R[V',\frac{s'V'+1}{r'}]$ (where we identified $U=\frac{s'V'+1}{r'}$). It is not that difficult to see that then there exist $a,b,c\in R$ 
such that $\alpha V'+\beta=aV'+b\frac{s'V'+1}{r'}+c$. This means that $\alpha=a+b\frac{s'}{r'}$, thus $\lambda \frac{\ttr}{r'}=a+b\frac{s'}{r'}$. 
This means that $\lambda\ttr=ar'+bs'$, and since $\lambda\in R^*$ this means $\ttr\in r'R+s'R$. Of course, the same method will also yield
$\tts\in r'R+s'R, r',s'\in \ttr R+\tts S$, hence the ideals $(\ttr,\tts)$ and $(r',s')$ are equal. The theorem is proved.
\end{proof}

\section{Conclusions and new examples}

Combining \ref{gelijk} and \ref{ongelijk} it is possible to construct a wider class of UFD counterexamples to generalized cancellation.
To give a new example, take $R$ a rigid ring from \cite{FM06}, like $R:=\C[X,Y,Z]/(X^2+Y^3+Z^7)$. (There are few rings 
known to be rigid! That's why we recycle this ring.) Now choose $r:=p(x), s=q(y),r':=\tilde{p}(x), s'=\tilde{q}(y)$ 
where $p,q, \tilde{p}, \tilde{q}$ are polynomials in one variable. 
Require that $p,\tilde{p}$ (resp. $q,\tilde{q}$) have the same zeroes (i.e. their radicals are the same), to make sure that they are stably isomorphic.
Possible choices are $p=x(x-1),q=y, \tilde{p}=x^2(x-1), \tilde{q}=y$, but also $p=x, q=y, \tilde{p}=2x, \tilde{q}=y$. 
In \cite{FM06} it is shown that an automorphism of $R$ sends $(x,y,z)$ to $(\lambda x, \mu y, \nu z)$ where $\lambda,\mu, \nu \in \C$. This can be used 
to show that there exists no automorphism sending $p$ to $\tilde{p}$ and $q$ to $\tilde{q}$ in general. 
In particular, the case $p=x(x-1),q=y, \tilde{p}=x^2(x-1), \tilde{q}=y$
gives a new counterexample to generalized cancellation. 

As mentioned before, it is not possible this way to find a counterexample to ``the'' cancellation problem ( If $A^{[1]}=\C^{[n]}$, then $A\cong \C^{[n-1]}$)
as $A_{r,s}$ can never be a polynomial ring. However, the reader may wonder if some of the choices made in section 2 can be improved upon.\\

{\bf Acknowledgements:} The author would like to thank both prof. dr. Finston and the anonymous referee for some useful suggestions and corrections.\\

\end{document}